\documentclass[12pt]{amsart}
\usepackage{amscd,amssymb}
\usepackage[plainpages,backref,urlcolor=blue]{hyperref}
\usepackage[all]{xy}

\theoremstyle{plain}
\newtheorem{thm}[subsection]{Theorem}
\newtheorem{lem}[subsection]{Lemma}
\newtheorem{prop}[subsection]{Proposition}
\newtheorem{cor}[subsection]{Corollary}

\theoremstyle{definition}
\newtheorem{rk}[subsection]{Remark}
\newtheorem{definition}[subsection]{Definition}
\newtheorem{ex}[subsection]{Example}
\newtheorem{conj}[subsection]{Conjecture}

\numberwithin{equation}{section}
\newcommand{\OO}{{\mathcal O}}

\newcommand{\ZZ}{{\mathcal Z}}

\newcommand{\I}{{\mathcal I}}
\newcommand{\J}{{\mathcal J}}
\newcommand{\F}{{\mathcal F}}

\newcommand{\CC}{{\mathcal C}}

\newcommand{\E}{{\mathcal E}}
\newcommand{\V}{{\mathcal V}}

\newcommand{\al}{{\alpha}}
\newcommand{\be}{{\beta}}

\newcommand{\C}{\mathbb{C}}
\newcommand{\T}{\mathbb{T}}
\newcommand{\PP}{\mathbb{P}}

\DeclareMathOperator{\Hom}{Hom}

\begin{document}

\title [Deformations of plane curves and Jacobian syzygies]
{Deformations of plane curves and Jacobian syzygies}

\author[Alexandru Dimca]{Alexandru Dimca}
\address{Universit\'e C\^ ote d'Azur, CNRS, LJAD, France }
\email{dimca@unice.fr}

\author[Gabriel Sticlaru]{Gabriel Sticlaru}
\address{Faculty of Mathematics and Informatics,
Ovidius University
Bd. Mamaia 124, 900527 Constanta,
Romania}
\email{gabrielsticlaru@yahoo.com }





\begin{abstract} 
We relate the equianalytic and the equisingular deformations of a reduced complex plane curve to the Jacobian syzygies of its defining equation.
Several examples and conjectures involving rational cuspidal curves are discussed.
\end{abstract}
 
\maketitle


\section{Introduction}
Let $C:f=0$ be a reduced plane curve in the complex projective plane $\PP^2$. The equianalytic and the equisingular deformations of a plane curve singularity $(C,p)$ are rather well understood, and the corresponding theory is briefly recalled in the second section following \cite{DiHa88,Wa74}. On the other hand, the  equianalytic and the equisingular deformations of the plane curve $C$, encoded in the (possibly non-reduced) analytic subspaces $\V_d^*(S_1,...,S_r)$ of the projective space $\PP^N$, are well understood when $C$ is nodal, see \cite{H, SeBook, Sev},
but much less in the general case, see \cite{DeSe, GK, GL96, GL01,GLS, GLSnb, KS,Lue,SeBook, Ta, Wa74A}. Here $d$ is the degree of the curve $C$, $S_1$,...,$S_r$ is the list of {\it all}  the singularities of $C$, $*=ea$ or $*=es$ denotes the equianalytic or the equisingular deformations, and $N=d(d+3)/2$, such that $\PP^N$ parametrizes all the degree $d$ plane curves. The main general results on these analytic subspaces
$\V_d^*(S_1,...,S_r)$ are recalled below in Theorem \ref{thm1}.

In the third section, we study the subspaces $\V_d^{ea}(S_1,...,S_r)$
by restating Theorem \ref{thm1} in terms of the Jacobian syzygies of $f$ and its Jacobian module $N(f)$, following a result by E. Sernesi, see
\cite[Corollary 2.2]{Se}. The main results here are Theorem \ref{thm3}, describing the tangent space $T_C\V_d^{ea}(S_1,...,S_r)$, and Theorem \ref{thm4}, giving a necessary and sufficient condition for $C$ to be unobstructed in terms of $mdr(f)$, the minimal degree of a Jacobian syzygy for $f$. 

In this paper we are particularly interested in the deformation theory of rational cuspidal plane curves. As it is conjectured that any such curve is either free or nearly free, see Conjecture \ref{c1}, we state a number of corollaries for free and nearly free curves, see for instance Corollary \ref{cor1} and Corollary \ref{cor3}, and we discuss such curves in a number of examples, see for instance Example \ref{ex2}.
In the fourth section, we recall first the classification of rational cuspidal plane curves with at least 3 cusps involving three infinite families $\F\ZZ_1(d,a)$, $\F\ZZ_2(k)$ and $\F\E(k)$ introduced by Flenner, Zaidenberg and Fenske in \cite{Fe99, FZ96,FZ00}. Concerning the curves in these three families, we conjecture that all of them are {\it free}, and that the corresponding invariant $mdr(f)$ decides to which family the curve $C:f=0$ belongs, see Conjecture \ref{conjMC} for details. This conjecture was verified for the first few curves in each of the three families, using their description in terms of equations or parametrizations. Section 4 is completed by a brief discussion of the Rigidity Conjecture by Flenner and Zaidenberg, which implies in particular that a rational cuspidal plane curve $C$ with at least 3 cusps
is strongly $es$-rigid, i.e. $\V_d^{es}(S_1,...,S_r)$ is smooth and coincide with the $G$-orbit of $C$, where $G=PGL(3,\C)$ is the automorphism group of $\PP^2$, acting in the obvious way on $\PP^N$.

In the final section we discuss the equisingular deformations of a reduced plane curve, the main result being Theorem \ref{thm10}.
The subtle structure of the analytic subspace $\V_d^*(S_1,...,S_r)$ is highlighted in Examples \ref{ex35} and  \ref{ex36}, where the rational unicuspidal curve $C:f=y^d+x^{d-1}z=0$ is deformed, for $d=5$, and
respectively $d=6$. The corresponding subspace $\V_5^{ea}(S_1)$ consists of two $G$-orbits, and it is smooth, while  (the support of) $\V_6^{ea}(S_1)$
consists of a single $G$-orbit, but $\V_6^{ea}(S_1)$ is not reduced
at all points. This latter fact is similar to classical examples of B. Segre of plane curves with many cusps, see for instance \cite[Corollary 2.3]{Ta} or \cite[Example 4.7.10]{SeBook}, and to Wahl's example of a curve of degree $d=104$ having $3636$ nodes and $900$ cusps, see \cite{Wa74A} and \cite[Examples 6.4 (6)]{GK}, but surprisingly in our example
the dimension of $\V_6^{ea}(S_1)$ is the expected dimension.
The corresponding equisingular deformations subspaces
$\V_5^{es}(S_1)$ and $\V_6^{es}(S_1)$ are unions of the orbit
$G \cdot C$ with families of $G$-orbits parametrized by weighted projective spaces of dimension 1 and 2 respectively. Each orbit in these families contains $G \cdot C$ in its closure, such that both analytic subspaces $\V_5^{es}(S_1)$ and $\V_6^{es}(S_1)$ turn out to be smooth of the expected dimension.

\bigskip

The authors thank Gert-Martin Greuel, Karol Palka and Edoardo Sernesi for useful discussions and remarks concerning this paper.

\section{Various prerequisites}

\subsection{Basic facts on deformations of isolated plane curve singularities} 

Consider an isolated curve singularity $(C,0):g=0$ at the origin of $\C^2$, defined by an analytic function germ $g \in \OO_2=\C\{u,v\}$.
Let $I^{ea}(C,0)$ denote the ideal in $\OO_2$ generated by $g$ and by the partial derivatives $g_u,g_v$ of $g$ with respect to $u$ and $v$.
Then the Tjurina algebra of $g$ is the quotient $T(g)=\OO_2/I^{ea}(C,0)$ and its dimension, denoted by $\tau(g)$ or $\tau(C,0)=\tau^{ea}(C,0)$, is called the Tjurina number of $g$ or of $(C,0)$. The ideal $I^{ea}(C,0)$ is called the equianalytic ideal of $g$.
Moreover, it is known that the germ
$(T(g),0)$ is the base space of the miniversal deformation of the singularity $(C,0):g=0$, see for details \cite{DiHa88,GL96,GL01,GLS,Wa74}.

There is a smooth germ $(ES(g),0) \subset (T(g),0)$ corresponding to the equisingular deformations of the singularity $(C,0)$, and the tangent space $T_0ES(g) \subset T(g)=T_0T(g)$ is given by a quotient 
$$I^{es}(C,0)/I^{ea}(C,0),$$
where $I^{es}(C,0)$ is the equisingular ideal of $g$ in $\OO_2$, see for details \cite{DiHa88, GL96,GL01,GLS,Wa74}. The following example is taken from \cite{Wa74}.

\begin{ex}
\label{ex1}
Consider the weighted homogeneous singularity $(C,0):u^p+v^q=0$ with $2 \leq p \leq q$. We will denote this singularity by $T_{p,q}$ following \cite{KS}.
Then the equianalytic ideal $I^{ea}(C,0)$ is generated by
$u^{p-1}$ and $v^{q-1}$, while the equisingular ideal  $I^{es}(C,0)$
is generated by $u^{p-1}, v^{q-1}$ and all the monomial $u^ a v^b$ with 
$aq+bp \geq pq$. This corresponds to the fact that a family of plane curve singularities over a reduced base is equisingular if and only if it is $\mu$-constant, and the Milnor number of an isolated weighted homogeneous singularity is not changed by adding higher order terms.
\end{ex}

\begin{rk}
\label{rk1}
Note that we have 
$$I^{es}(C,0)=I^{ea}(C,0)$$
if  $(C,0)$ is a simple singularity, i.e. it is a singularity of type $A_n$, $D_n$, $E_6$, $E_7$ or $E_8$ in Arnold's classification.
\end{rk}

One defines the {\it equisingular Tjuriana number} of the singularity $(C,0)$ by the equality
\begin{equation}
\label{eq0}
\tau^{es}(C,0)=\dim \OO_2/I^{es}(C,0),
\end{equation}
and the {\it equisingular modality} of the singularity $(C,0)$ by the equality

\begin{equation}
\label{eq01}
m^{es}(C,0)=\dim I^{es}(C,0)/I^{ea}(C,0)=\tau(C,0)-\tau^{es}(C,0).
\end{equation}
The modality $mod(g)$ of the germ $g$ with respect the the right-equivalence is also related to these invariants by the formula 
\begin{equation}
\label{eq02}
mod(f)=\mu(C,0)-\tau^{es}(C,0),
\end{equation}
see \cite[Lemma 1.7]{GL96}.
\subsection{Basic facts on deformations of reduced plane curves} 

Let $S=\C[x,y,z]$ be the polynomial ring in three variables $x,y,z$ with complex coefficients and let $\PP(S_d)$ be the projective space parametrizing the degree $d$ curves in $\PP^2$. Let $N=\dim \PP(S_d)=d(d+3)/2$.
Fix a list $S_1,S_2,...,S_r$ of isomorphism classes of isolated plane curve singularities. We consider two settings, namely
the equianalytic setting, denoted with ${}^{ea}$, and the equisingular setting, denoted with ${}^{es}$. We denote by $V_d^*(S_1,...,S_r)$ the subset in $\PP(S_d)$ corresponding to the reduced curves $C \subset \PP^2$, having precisely $r$ singularities which are required to be
\begin{enumerate}

\item isomorphic to the singularities $S_1,...,S_r$ when $*=ea$, and

\item equisingular to the singularities $S_1,...,S_r$ when $*=es$.
\end{enumerate}
By definition we have 
\begin{equation}
\label{eq1}
V_d^{ea}(S_1,...,S_r)\subset V_d^{es}(S_1,...,S_r),
\end{equation}
with equality when all the singularities $S_j$'s are simple. 

Each of the two sets $V_d^*(S_1,...,S_r)$ has in fact the structure of a natural (possibly non reduced) complex analytic subspace $\V_d^*(S_1,...,S_r)$ of the projective space $\PP^N$, such that the restriction of the universal family over $\PP^N$
to the analytic subspace $\V_d^*(S_1,...,S_r)$ has the expected functorial properties. The first rigorous proof of the existence of this analytic subspace structure was given by J. Wahl in \cite{Wa74A}, for curves having only nodes $A_1$ and simple cusps $A_2$ as singularities.
The general construction for the analytic subspace structure of $\V_d^{ea}(S_1,...,S_r)$ was settled by G.-M. Greuel and U. Karras in \cite{GK},
see Theorem 1.3 and Theorem 2.2. This proof was then adapted for 
the analytic subspace structure of $\V_d^{es}(S_1,...,S_r)$ by G.-M. Greuel and C. Lossen in \cite{GL96}. Another, more direct approach to the existence of the analytic subspace structure of $\V_d^{*}(S_1,...,S_r)$ is given in the forthcoming book \cite{GLSnb}, see Theorems II. 2.32 and II.2.36, where the authors also show that the reduction of these subspaces
are quasi-projective subsets of $\PP^N$. It is very likely that the
analytic subspace structures $\V_d^{*}(S_1,...,S_r)$ themselves are (algebraic) subschemes of $\PP^N$, but this fact is proved only in the case of curves with nodes and cusps as singularities, see \cite[Theorem 3.3.5]{Wa74A}.

Consider now a reduced curve $C:f=0$ in $\PP^2$, having $r$ singular points located at the points $p_j \in \PP^2$, for $j=1,...,r$. To such a curve we associate two 0-dimensional subschemes of $\PP^2$, namely
$Z^*(C)$ for $*=ea,es$, such that the support of the schemes $Z^*(C)$
is the singular points of $C$ and one has by definition
$$\OO_{Z^*(C),p_j}=\OO_{\PP^2,p_j}/I^*(C,p_j).$$
By our discussion in the previous subsection, it follows that
$$\deg Z^*(C)= \sum_{j=1,r} \tau^*(C,p_j).$$

Let $\I(Z^*(C))$ be the ideal sheaf in $\OO_{\PP^2}$ which defines the scheme $Z^*(C)$. Note that $Z^*(C)$ can also be regarded as a subscheme of the curve $C$, and  then it is defined by the ideal sheaf
$$\I(Z^*(C)|C)=\I(Z^*(C)) \otimes \OO_C \subset \OO_C.$$
As noted in \cite[Remark 2.4]{GL01}, one has an exact sequence
$$0 \to \I_C(d) \to \I(Z^*(C))(d) \to \I(Z^*(C)|C)(d) \to 0,$$
where $\I_C \subset \OO_{\PP^2}$ denotes the sheaf ideal of $C$,
and one has $\OO_{\PP^2} \simeq  \I_C(d)$. It follows that
$$H^0(\I(Z^*(C)|C)(d))=H^0(\I(Z^*(C))(d))/H^0(\OO_{\PP^2}).$$
These sheaf ideals are important in view of the following result.
\begin{thm}
\label{thm1}
Let $C$ be a reduced plane curves such that $C \in V_d^*(S_1,...,S_r)$, with $*=ea,es$. Then the following hold.
\begin{enumerate}

\item $H^0(\I(Z^*(C)|C)(d))$ is isomorphic to the Zariski tangent space of
the analytic subspace $ \V_d^*(S_1,...,S_r)$ at the point $C$.

\item $H^1(\I(Z^*(C)|C)(d))=0$ if and only if the analytic subspace $ \V_d^*(S_1,...,S_r)$ is smooth at the point $C$, of the expected dimension, namely
$$N-\deg Z^*(C)= d(d+3)/2- \sum_{j=1,r} \tau^*(C,p_j).$$

\end{enumerate}

\end{thm}

The  claim (1) above is not very difficult and can be found in \cite[Proposition 4.19]{DiHa88}. The second claim (2) is more difficult, and we refer to \cite{GL01} and to the forthcoming book \cite{GLSnb}, see Theorems II. 2.38 and II.2.40 for the complete proofs.

\subsection{Jacobian ideal, Jacobian module, and free and nearly free curves}

We denote by $J_f$ the Jacobian ideal of $f\in S_d$, i.e. the homogeneous ideal in $S$ spanned by the partial derivatives $f_x,f_y,f_z$, and  by $M(f)=S/J_f$ the corresponding graded quotient ring, called the Jacobian (or Milnor) algebra of $f$.
 Let $I_f$ denote the saturation of the ideal $J_f$ with respect to the maximal ideal ${\bf m}=(x,y,z)$ in $S$ and consider the local cohomology group, usually called the Jacobian module of $f$, 
 $$N(f)=I_f/J_f=H^0_{\bf m}(M(f)).$$
 The graded $S$-module $AR(f)=AR(C) 
\subset S^{3}$ of 
\textit{all Jacobian relations} of $f$ is defined by 
\begin{equation} \label{eqAR} 
AR(f)_k:=\{(a,b,c) \in S^{3}_k \mid a f_x+b f_y+c f_z=0\}.
\end{equation}
Its sheafification $E_C:=
\widetilde{AR(f)}$ is a rank two vector bundle on $\PP^2$, see 
\cite{DS14,Se} for details. More precisely, one has $E_C=T\langle C \rangle (-1)$,
where $T\langle C \rangle $ is the sheaf of logarithmic vector fields along $C$ as considered for instance in \cite{DS14,Se}. We set 
$$ar(f)_m= \dim AR(f)_m=\dim H^0(\PP^2,E_C(m)) \text{ and } n(f)_m= \dim N(f)_m$$
 for any integer $m$.
The minimal degree of a Jacobian relation for the polynomial $f$ is the integer $mdr(f)$
defined to be the smallest integer $m\geq 0$ such that $ar(f)_m >0$.
One clearly has $mdr(f)=0$ if and only if, up-to a linear change of coordinates in $G$, $f$ is independent of $z$, i.e. the curve $C$ consists of $d$ concurrent lines. 

\medskip

\noindent {\it We assume from now on in this paper that} $mdr(f)>0$.

\medskip

\noindent Moreover, the irreducible curves $C$ as well as the line arrangements with $mdr(f)=1$ are classified, see \cite{duPCTC2}. Hence the interesting cases are the curves with $ar(f)_1=0$.

It was shown in \cite[Corollary 4.3]{DPop} that the graded $S$-module  $N(f)$ satisfies a Lefschetz type property with respect to multiplication by generic linear forms. This implies in particular the inequalities
$$0 \leq n(f)_0 \leq n(f)_1 \leq ...\leq n(f)_{[T/2]} \geq n(f)_{[T/2]+1} \geq ...\geq n(f)_T \geq 0,$$
where $T=3d-6$. We set as in \cite{Drcc}
$$\nu(C)=\max _j \{n(f)_j\},$$
and introduce a new invariant for $C$, namely
$$\sigma(C)=\min \{j   : n(f)_j \ne 0\}.$$
The self duality of the graded $S$-module $N(f)$, see \cite{Se, SW}, implies that $n(f)_s \ne 0$ exactly for $s=\sigma(C),..., T-\sigma(C)$. 

Recall that $C$ is a free curve if $J_f=I_f$, or equivalently $\nu(C)=0$, see \cite{Dmax, DStFD,ST}. Note that for a free curve $C:f=0$ of degree $d$, one has the exponents $(d_1,d_2)$ where $d_1=mdr(f)$ and $d_2=d-1-d_1$, as well as the formula
\begin{equation}\label{eqTAU1}
\tau(C)=(d-1)^2-d_1d_2.
\end{equation}
Similarly, $C$ is a nearly free curve if $\nu(C)=1$, see \cite{B+,Dmax, DStRIMS}. Note that one has
$\sigma(C)=d+mdr(f)-3$ for a nearly free curve by \cite[Corollary 2.17]{DStRIMS}.
Moreover, for a nearly free curve $C:f=0$ of degree $d$, one has the exponents $(d_1,d_2)$ where $d_1=mdr(f)$ and $d_2=d-d_1$, and the formula
\begin{equation}\label{eqTAU2}
\tau(C)=(d-1)^2-d_1(d_2-1)-1.
\end{equation}
Our interest in the free and nearly free curves comes from the following.
\begin{conj}
\label{c1}
A reduced plane curve $C:f=0$ which is rational cuspidal is either free, or nearly free.
\end{conj}
This conjecture is known to hold when the degree $d$ of $C$ is even, or when $d \leq 33$, as well as in many other cases, see \cite{ Drcc, DStRIMS,DStMos}.

\section{Jacobian syzygies, Jacobian module and equianalytic deformations} 

The following result is obvious, but plays a key role in the sequal.
\begin{lem}
\label{lem1}
The ideal sheaf $\J_f$ associated to the homogeneous ideal $J_f$, coincides with the ideal sheaf $\I_f$ associated to the homogeneous ideal $I_f$, and with the ideal sheaf $\I(Z^{ea}(C))$. This ideal sheaf defines the scheme structure of the singular locus of the plane curve $C$.
\end{lem}

Denote by $\J_f|C=\J_f \otimes \OO_C$ the resctriction of the ideal sheaf $\J_f$ to the curve $C$ and note that in view of Lemma \ref{lem1} we have
$$\J_f|C=\I(Z^{ea}(C)|C).$$
Hence, Theorem \ref{thm1} can be restated in the case $*=ea$ as follows.
\begin{thm}
\label{thm2}
Let $C$ be a reduced plane curves such that $C \in V_d^{ea}(S_1,...,S_r)$. Then the following hold.
\begin{enumerate}

\item $H^0((\J_f|C)(d))$ is isomorphic to $T_C \V_d^{ea}(S_1,...,S_r)$,
the Zariski tangent space of the analytic subspace
$ V_d^{ea}(S_1,...,S_r)$ at the point $C$.

\item $H^1(\J_f|C(d))=0$ if and only if the analytic subspace $ \V_d^{ea}(S_1,...,S_r)$ is smooth at the point $C$, of the expected dimension 
$d(d+3)/2-\tau(C),$
where $\tau(C)=\sum_{j=1,r} \tau(C,p_j)$ is the total Tjurina number of the curve $C$.
\end{enumerate}

\end{thm}
In this section we look at the dimensions $h^0((\J_f|C)(d))=\dim H^0((\J_f|C)(d))$ and $h^1(\J_f|C(d))=\dim H^1(\J_f|C(d))$, using results by Sernesi in \cite{Se} and our study of the
Jacobian syzygies and  Jacobian module $N(f)$ in \cite{DS14}.

We recall the following diagram from \cite{Se}, with exact rows and columns.
\begin{equation}
\label{dia1}
	\xymatrix{&&0&0 \\
	0\ar[r]& T\langle C \rangle \ar[r] & T_{\PP^2}\ar[u] \ar[r]^-\eta& (\J_f|C)(d)\ar[u] \ar[r] & 0 \\
	0\ar[r]&E_C(1)\ar[u]^-{\cong}\ar[r]&\mathcal{O}_{\PP^2}(1)^{3} \ar[u]\ar[r]^-{\partial f} & \mathcal{J}_f(d) \ar[r]\ar[u] &0 \\
	&&\mathcal{O}_{\PP^2}\ar[u]\ar[r]^-f& \mathcal{I}_C(d)\ar[u]\ar[r]&0 \\
	&&0\ar[u]&0\ar[u]}
\end{equation}
Here the morphism $\partial f$ is induced by the morphism of graded $S$-modules
$$ S(1)^3 \to J_f(d), \  \  (a,b,c) \mapsto af_x+bf_y+cf_z,$$
and $\eta$ is induced by $\partial f$. The middle column is the usual free resolution of the tangent bundle $T_{\PP^2}$ to $\PP^2$, and implies
$$h^0(T_{\PP^2})=3h^0(\mathcal{O}_{\PP^2}(1))-h^0(\mathcal{O}_{\PP^2})=3 \cdot 3-1=8=\dim PGL(3,\C).$$
The first exact row yields the following exact sequence
$$0 \to AR(f)_1 \to H^0(T_{\PP^2}) \to H^0(\J_f|C)(d)) \to H^1(T\langle C \rangle) \to 0.$$
Indeed, one has $H^0(T\langle C \rangle)=H^0(E_C(1))=AR(f)_1$
and $H^1(T_{\PP^2})=0$ using the middle column in the above diagram.
Note that one has the following facts.

\begin{enumerate}

\item The automorphism group of $\PP^2$ is $G=PGL(3,\C)$, and this group acts in an obvious way on the projective space $\PP(S_d)$, such that all the subsets $ V_d^{*}(S_1,...,S_r)$ are unions of $G$-orbits.
In particular, one has the inclusions
$$T_C (G \cdot C) \subset T_C V_d^{ea}(S_1,...,S_r) \subset T_C \V_d^{ea}(S_1,...,S_r).$$
\item The tangent space $T_eG$ to $G$ at the identity element $e \in G$
can be identified to $H^0(T_{\PP^2})$;

\item Let $H \subset G$ be the subgroup of elements fixing the polynomial $f \in \PP(S_d)$. 
Then the tangent space $T_eH$ to $H$ at the identity element $e \in G$
can be identified to $H^0(T\langle C \rangle)=H^0(E_C(1))=AR(f)_1$, see for instance \cite{duPCTC2}.

\item Hence the cokernel of the inclusion $AR(f)_1 \to H^0(T_{\PP^2})$ can be identified to the tangent space $T_C (G \cdot C)$ to the orbit of $C$, at the point $C$.

\end{enumerate}
It follows from \cite[Proposition 2.1]{Se} that one has $H^1(T\langle C \rangle)=N(f)_d$. Therefore we get the following result, which is our reformulation of \cite[Corollary 2.2]{Se}.

\begin{thm}
\label{thm3}
Let $C$ be a reduced plane curves such that $C \in V_d^{ea}(S_1,...,S_r)$. Then  the Zariski tangent space
 $T_C \V_d^{ea}(S_1,...,S_r)$ of the analytic subspace
$ \V_d^{ea}(S_1,...,S_r)$ at the point $C$ sits in the following exact sequence
$$0 \to T_C (G \cdot C) \to T_C \V_d^{ea}(S_1,...,S_r) \to N(f)_d \to 0.$$
In particular, one has
$\dim T_C \V_d^{ea}(S_1,...,S_r)=8-ar(f)_1+n(f)_d$.
\end{thm}

\begin{definition}
\label{def1}
We say that the reduced plane curve $C$ is $*$-projectively rigid
if $T_C (G \cdot C) = T_C \V_d^{*}(S_1,...,S_r)$.
Note that this is equivalent to the equality of analytic space germs
$(G \cdot C,C)=(V_d^{*}(S_1,...,S_r),C) =(\V_d^{*}(S_1,...,S_r),C) .$
We say that the reduced plane curve $C$ is strongly $*$-projectively rigid, if the analytic subspace $\V_d^{*}(S_1,...,S_r)$ containing $C$ is smooth and coincides with the orbit $G \cdot C$.
\end{definition}
In particular, if $C$ is $*$-projectively rigid, it follows that the analytic subspace
$\V_d^{*}(S_1,...,S_r)$ is reduced and smooth at $C$.
Clearly, if a curve $C$ is (strongly) $es$-rigid, then $C$ is also (strongly) $ea$-rigid. We have also the following result, with an obvious proof.

\begin{lem}
\label{lem2} A reduced plane curve $C$ is (strongly) $*$-projectively rigid if and only if any curve $C'$ in the $G$-orbit $G \cdot C$ is (strongly) $*$-projectively rigid. If the reduced plane curve $C\in V_d^{*}(S_1,...,S_r)$ is $*$-projectively rigid, then the 
$G$-orbit $G \cdot C$ is an irreducible and connected component of the analytic subspace $\V_d^{*}(S_1,...,S_r)$, consisting only of smooth points of this analytic subspace.
\end{lem}
\begin{cor}
\label{cor1}
A reduced curve $C:f=0$ is $ea$-rigid if and only if $n(f)_d=0$.
In particular, one has the following.
\begin{enumerate}

\item Any free curve is $ea$-rigid.

\item A nearly free curve $C:f=0$ is $ea$-rigid if and only if one of the following cases occurs: 

\subitem (i) $mdr(f) \geq 4$, or

\subitem (ii) $mdr(f)=2$ and $d=4$, or

\subitem (iii) $mdr(f)=1$ and $d=2,3$.

\end{enumerate}

\end{cor}
The  claim (2) above  corrects the final claim in \cite[Corollary 2.17]{DStRIMS}.

\proof
To prove the second claim, note that $n(f)_d=0$ if either
$d<\sigma(C)=d+mdr(f)-3$ or $d>T-\sigma(C)$. The first inequality is equivalent to $mdr(f) \geq 4$. The second inequality is equivalent to
$mdr(f)>d-3$. For a nearly free curve one has $mdr(f)\leq d/2$, see \cite{DStRIMS}, and this gives the subcases (ii) and (iii) above. The following examples shows that the sub cases (ii) and (iii) above really occur.

\endproof

\begin{ex}
\label{ex2} In this example we list the free and nearly free curves $C:f=0$ of degree $d$ such that $2 \leq d \leq 4$ and $mdr(f)>0$.
For details we refer to \cite{KS, Moe}. All of them are strongly $ea$-rigid except for the curve $C_4'$ discussed below.
\begin{enumerate}

\item For $d=2$, the smooth conic $C:f=x^2+y^2+z^2=0$ is nearly free,
with $ar(f)_1=3$ and $mdr(f)=1$. The scheme $\V_2^{*}(\emptyset)$ coincides with the orbit $G \cdot C$, and it is an Zariski open subset in 
$\PP^5=\PP(S_2)$. Moreover, $n(f)_k=0$ for $k \ne 0$, so the formula for $\dim \V_2^{*}(\emptyset)$ in Theorem \ref{thm3} holds.

\item For $d=3$, the cuspidal cubic $C:f=x^2y+z^3=0$ is nearly free,
with $ar(f)_1=1$ and $mdr(f)=1$. The quasi-projective variety $V_3^{*}(A_2)$ coincides with the orbit $G \cdot C$, and it has codimension 2 in
$\PP^9=\PP(S_3)$. Moreover, $n(f)_3=0$, so the formula for $\dim \V_3^{*}(A_2)$ in Theorem \ref{thm3} shows that the analytic subspace $\V_3^{*}(A_2)$ is smooth. The triangle $C':f=xyz=0$ is free, with $ar(f)_1=2$ and $mdr(f)=1$. The analytic subspace $\V_3^{*}(3A_1)$ is smooth, coincides with the orbit $G \cdot C'$, and has dimension 6.

\item For $d=4$, there are 4 lists of singularities which may occur on a nearly free irreducible quartic, namely $3A_2$, $A_4A_2$, $A_6$ and $E_6$. The analytic subspace $\V_4^{*}(3A_2)$, $\V_4^{*}(A_4A_2)$ and $\V_4^{*}(A_6)$ are smooth and coincide with the corresponding $G$-orbits, while $\V_4^{*}(E_6)$ is still smooth, but the  union of two $G$-orbits. Consider the two quartics
$$C_4:f=y^4-xz^3-y^3z=0 \text{ and } C_4':f'=y^4-xz^3=0.$$ 
For all these quartics, except for $C_4'$, one has $mdr(f)=2$, which implies $ar(f)_1=n(f)_4=0$, and hence all these $G$-orbits are 8-dimensional in $\PP^{14}=\PP(S_4)$. For $C_4'$ one has $ar(f)_1=n(f)_4=1$,
and hence $C_4'$ is not $ea$-rigid, since
$$\dim T_{C_4'}(G \cdot C_4')=7 <8=\dim T_{C_4'}\V_4^{*}(E_6).$$
Note also that the orbit $G \cdot C_4'$ is contained in the closure of the orbit $G \cdot C_4$.
It follows that the variety $V_4^{*}(E_6)$ is the union of these two orbits,
and hence the corresponding analytic subspace  $\V_4^{*}(E_6)$ is smooth and irreducible of dimension 8 at any point.

The line arrangement $C':f=xyz(x+y+z)=0$ is also nearly free, the orbit $G\cdot C'$ is 8-dimensional as above, and coincides with $\V_4^{*}(6A_1)$. The line arrangement $C'':f=xyz(x+y)=0$ is free with $mdr(f)=ar(f)_1=1$, the orbit $G\cdot C''$ is 7-dimensional, and coincides with $\V_4^{*}(D_43A_1)$. 
Note that $E_6=T_{3,4}$ and $D_4=T_{3,3}$.

\end{enumerate}

\end{ex}

Hence we have seen that the free and some of the nearly free curves give rise to smooth components of the analytic subspace $\V_d^{ea}(S_1,...,S_r)$. Now we investigate when such components have the expected dimension, given by Theorem \ref{thm2} (2).

\begin{thm}
\label{thm4}
Let $C$ be a reduced plane curves such that $C \in V_d^{ea}(S_1,...,S_r)$. Then the analytic subspace
$ \V_d^{ea}(S_1,...,S_r)$ is smooth at the point $C$ of the expected dimension 
$$\frac{d(d+3)}{2}-\tau(C)$$
 if and only if $mdr(f) \geq d-4$. 
\end{thm}

\proof Using the first row in the commutative diagram \eqref{dia1}, we see that there is an isomorphism $H^1(\J_f|C(d))=H^2(T\langle C \rangle)$.
Then, using Serre duality and the identity
$$\Omega(\log C)=\Hom (T\langle C \rangle, \OO_{\PP^2})=T\langle C \rangle(d-3)$$
as in \cite[Lemma 4.1]{Se}, we get
$$h^2(T\langle C \rangle)=h^0(T\langle C \rangle(d-6)=h^0(E_C(d-5))=ar(f)_{d-5}.$$

\endproof

\begin{cor}
\label{cor2}
Let $C$ be a reduced plane curves such that $C \in V_d^{ea}(S_1,...,S_r)$.
Then the analytic subspace
$ \V_d^{ea}(S_1,...,S_r)$ is smooth at the point $C$ of the expected dimension 
$$\frac{d(d+3)}{2}-\tau(C)$$
 if the curve $C$ satisfies one of the following conditions.

\begin{enumerate}

\item The degree $d$ of $C$ is at most 5.

\item The degree $d$ of $C$ is at most 6 and $C$ is free and irreducible.

\item The curve $C$ is a nodal curve.

\end{enumerate}

\end{cor}
The third claim is known, see \cite{H,Sev}, but our approach is perhaps new.
\proof
The first claim follows from our assumption that we consider in this paper only curves with $mdr(f)>0$. The second claim follows from the fact that an irreducible free curve $C:f=0$ satisfy $mdr(f) \geq 2$, see \cite[Theorem 2.8]{DStFD}.
The third claim follows from the fact that a nodal curve $C:f=0$ satisfy $mdr(f) \geq d-2$, see \cite[Theorem 4.1]{DStEdin}.
\endproof

\begin{cor}
\label{cor3}
Let $C$ be a reduced plane curves such that the analytic subspace
$ \V_d^{ea}(S_1,...,S_r)$ is smooth at the point $C$ of the expected dimension 
$$\frac{d(d+3)}{2}-\tau(C).$$
Then one has the following.

\begin{enumerate}

\item If the curve $C$ is free, then degree $d$ of $C$ is at most 7.

\item If the curve $C$ is nearly free, then degree $d$ of $C$ is at most 8.
 \end{enumerate}

\end{cor}

\proof It is enough to recall that for a free (resp. nearly free) curve $C:f=0$ one has
$mdr(f)<d/2$ (resp. $mdr(f) \leq d/2$).
\endproof

For a rational cuspidal plane curve $C$ which satisfies  Conjecture \ref{c1}, we can test whether $C$ is $ea$-rigid by using Corollary \ref{cor1},
and test the fact that $C$ is $ea$-unobstructed by using Theorem
\ref{thm4} and Corollary \ref{cor3}.

The next example shows that for curves which are not (nearly) free, the difference between the orbit $G \cdot C$ and the scheme $\V_d^{ea}(S_1,...,S_r)$ can be quite large.

\begin{ex}
\label{ex20} In this example we revisit some examples from \cite{GL96, Lue}
and treat them from our point of view.
\begin{enumerate}

\item The curve $C:f=y(x+2y+z)(x-2y-z)(x^4-x^2z^2+y^2z^2+y^3z)=0$
is considered in \cite[Example 5.3]{GL96}. This curve has three triple points of type $D_4$ and seven nodes $A_1$, hence $\tau(C)=3\cdot 4+7=19$ and $\V_7^{ea}(3D_4,7A_1)=\V_7^{es}(3D4,7A_1)$, since only simple singularities are involved. Greuel and Lossen show that the scheme $\V_7^{ea}(3D_4,7A_1)$ is smooth at $C$ of dimension 16.
This result follows also from  Theorem \ref{thm4}, since one has
$d=7$ and $mdr(f)=5\geq d-4=3$.

\item The curve $C:f=x^9+z(xz^3+y^4)^2=0$
is considered in \cite[Example 5.6 (c)]{GL96}. This curve is irreducible, has a unique singularity, which is of type $A_{35}$ and was used by Luengo in \cite{Lue} to construct the first singular variety $ V_d^{ea}(S_1,...,S_r)$. It is known that for this curve one has $h^0((\J_f|C)(d))=20$ and $h^1((\J_f|C)(d))=1$, see  \cite[Examples 6.4 (6)]{GK}.
We can recover these dimensions as follows. For this polynomial $f$, one has, by direct computation using SINGULAR, $d=9$, $n(f)_9=12$
and $mdr(f)=4$. By Theorem \ref{thm3} we have
$$h^0((\J_f|C)(d))=\dim T_C \V_9^{ea}(A_{35})=8+12=20.$$
From the proof of Theorem \ref{thm4}, we have that
$$h^1((\J_f|C)(d))=ar(f)_{4}=1,$$
since  in a minimal set of generators for the graded $S$-module $AR(f)$ there is only one of degree 4 and three more generators, all having degree 8.

\item The curve $C:f=x^9+x^8z+z(xz^3+y^4)^2=0$
is considered in \cite[Example 5.6 (c)]{GL96}. This curve is irreducible, has a unique singularity, which is of type $A_{31}$ and is a deformation of  Luengo's curve discussed above.  Greuel and Lossen show that the analytic subspace $\V_9^{ea}(A_{31})$ is smooth at $C$ of dimension 23.
This result follows also from  Theorem \ref{thm4}, since one has
$d=9$ and $mdr(f)=5\geq d-4=5$.

\end{enumerate}

\end{ex}

\section{On the rational plane curves with at least 3 cusps}
It is well known that a rational plane curve with at least 3 cusps has degree $d \geq 4$. Moreover, one has the following, see \cite{Moe,Nam}.

\begin{prop}\label{prop0C}
There is only one rational plane curve with  3 cusps of degree $d= 4$ up to projective equivalence, namely the tricuspidal quartic
$$\CC_4:[(2),(2),(2)]=[3A_2],  \  \  (s^3t-\frac{1}{2}s^4:s^2t^2:t^4-2st^3).$$
Up to projective equivalence, there are three rational cuspidal quintics
having at least 3 cusps. They have the following cuspidal configurations and parametrizations.

\begin{enumerate}

\item $\CC_5: [(3),(2_2),(2)]=[E_6,A_4,A_2]$, $ (s^4t-\frac{1}{2}s^5:s^3t^2:-\frac{3}{2}st^4+t^5)$.

\item $\CC_5':[(2_2),(2_2),(2_2)]=[3A_4]$, $(s^4t-s^5:s^2t^3-\frac{5}{32}s^5:-\frac{125}{128}s^5-\frac{25}{16}s^3t^2-5st^4+t^5)$.

\item $\CC_5'': [(2_3),(2),(2),(2)]=[A_6,3A_2]$, $(s^4t:s^2t^3-s^5:t^5+2s^3t^2)$.

 \end{enumerate}

\end{prop}
Using our algorithm to decide the freeness of a rational curve given by a parametrization described in \cite{BDS}, we get the following.
\begin{cor}\label{cor0C}
The curve $\CC_4$ is nearly free with exponents $(2,2)$.
The curves $\CC_5$, $\CC_5'$ and $\CC_5''$ are free with exponents $(2,2)$. Moreover, $\CC_4$, $\CC_5$ and $\CC_5''$ have the following equations
$$\CC_4: f(x,y,z)=x^2y^2+y^2z^2+x^2z^2-2xyz(x+y+z)=0,$$
$$\CC_5:f=9xy^4-4y^5-24x^2y^2z+48xy^3z-16y^4z+16x^3z^2=0,$$
$$\CC_5'':f=-27x^5+2x^2y^3-18x^3yz+y^4z-2xy^2z^2+x^2z^3=0.$$

\end{cor}
The equation for $\CC_4$ given above does not correspond
to the given parametrization. The equation for $\CC_5'$ obtained from the above parametrization is too complicated to list.

The rational cuspidal plane curves with 3 cusps of degree $d\geq 6$ are believed to by described by the following three infinite series, \cite{Pi}. 
The first series, call it $\F\ZZ_1(d,a)$, is described in \cite{FZ96}, see Theorem 3.5 and Proposition 3.9.
\begin{thm}
\label{thm1C}
Let $C$ be a rational cuspidal curve of type $(d,d-2)$ with at least three cusps. Then  $C$ has exactly three cusps  and there exists a unique pair of integers $a,b$, $a\geq b \geq 1$ with $a+b=d-2$ such that the multiplicity sequences of the cusps are $[(d-2),(2_a),(2_b)]$. Moreover,
up to projective equivalence the equation of $C=C_{d,a}$ can be written in affine coordinates $(x,y)$ as
$$f(x,y)= \frac{x^{2a+1}y^{2b+1}-((x-y)^{d-2}-xyg(x,y))^2}{(x-y)^{d-2}},$$
where $d \geq 4$,  $g(x,y)=y^{d-3}h(x/y)$ and 
$$h(t)= \sum_{k=0,d-3}\frac{a_k}{k!}(t-1)^k,$$
with $a_0=1$, $a_1=a-\frac{1}{2}$ and $a_k=a_1(a_1-1) \cdots (a_1-k+1)$ for $k>1$.
\end{thm}

\begin{conj}\label{conj1C}
The curves $C_{d,a} \in \F\ZZ_1(d,a)$  are free divisors and $mdr(f)=2$ for all possible values of the pair $(d,a)$, when $d\geq 5$.
In particular, $\tau(C_{d,a})=d^2-4d+7$ for any pair $(d,a)$ with $d\geq 5$.
\end{conj}
For $d=4$ the curve $C_{4,2}$ is the quartic with 3 cusps $A_2$ from Example \ref{ex2} (3), and it is nearly free with exponents $(d_1,d_2)=(2,2)$. This conjecture was verified by direct computation for all pairs $(d,a)$ with $5\leq d \leq 18$.

The  second series, call it $\F\ZZ_2(k)$, is described in \cite{FZ00}, see Theorem 1.1.
\begin{thm}
\label{thm2C}
Let $C$ be a rational cuspidal curve of type $(d,d-3)$ with at least three cusps and $d \geq 6$. Then $d=2k-1$ where $k\geq 4$, $C$ has exactly three cusps with multiplicity sequences  $[(2(k-2),2_{k-2}),(3_{k-2}),(2)]$. Moreover,
up to projective equivalence, this curve $C=C_k$  is given by a parametrization
$$(s:t) \mapsto (s^{2k-4}t^3:s^{2k-4}(s-t)^2(2s+t):t^3(s-t)^2q_k(s,t)),$$
where $q_k \in \C[s,t]$ is the homogeneous polynomial of degree $2k-6$ defined below.
\end{thm}
The polynomial $q_4$ is given in \cite{FZ00}, see equation \eqref{eqk2} below.
The polynomials $q=q_k$, for $k>4$,  are constructed as follows. We start with a polynomial $f \in \C[x,y]$, homogeneous of degree $k-5$. We make the substitution $x=t^3$ and $y=t^3-3t+2$, and we define a polynomial in $t$ by the formula
$$A(t)=f(t^3,t^3-3t+2)+t^{3(k-4)}g(t),$$
where $g(t)$ is a polynomial in $t$ of degree $2k-7$. The polynomials $f$ and $g$ are uniquely determined by the following conditions.
\begin{enumerate}

\item The $m$-th derivative of $A$ vanish at $t=-2$ for $0 \leq m \leq k-5$.

\item The $n$-th derivative of $A$ vanish at $t=1$ for $0 \leq n \leq 2k-9$.

\item The polynomial in $t$ defined by $$h(t)=\frac{A(t)}{(t^3-3t+2)^{k-4}}$$ 
satisfies $h(1)=-1$ and $h'(1)=3$.

 \end{enumerate}
 Finally, we define $q$ by asking that 
 $$q(1,t)=\frac{t^3h(t)+1}{(t-1)^2}.$$
Here are the first of these polynomials, $q_4$ being given in \cite{FZ00} and $q_5,q_6$ being computing using this approach.
\begin{equation}\label{eqk2}
q_4(s,t)=s^2+2st+3t^2.
\end{equation}

\begin{equation}\label{eqk3}
q_5(s,t)=s^4+2s^3t+3s^2t^2+\frac{36}{11}st^3+\frac{27}{11}t^4.
\end{equation}

\begin{equation}\label{eqk4}
q_6(s,t)=s^6+2s^5t+3s^4t^2+\frac{612}{169}s^3t^3+\frac{621}{169}s^2t^4+\frac{468}{169}st^5+\frac{243}{169}t^6.
\end{equation}

\begin{equation}\label{eqk5}
q_7(s,t)=s^8+2s^7t+3s^6t^2+\frac{3780}{1009}s^5t^3+\frac{4149}{1009}s^4t^4+\frac{3942}{1009}s^3t^5+\frac{3159}{1009}s^2t^6+\frac{1944}{1009}st^7+\frac{729}{1009}t^8.
\end{equation}

\begin{conj}\label{conj2C}
The curves $C_{k}\in \F\ZZ_2(k)$  are free divisors and $mdr(f)=k-1$ for all  $k \geq 4$.
In particular, $\tau(C_d)=3(k-1)^2$  for any $k\geq 4$.
\end{conj}
This conjecture was verified by direct computation for $4 \leq k \leq 7$ using the above parametrizations.

The  third series, call it $\F\E(k)$, is described in \cite{Fe99}, see Theorem 1.1. For the definition of $T_V(\langle D \rangle)$, see the final section.
\begin{thm}
\label{thm3C}
Let $C$ be a rational cuspidal curve of type $(d,d-4)$ with at least three cusps, $d \geq 10$ and satisfying $\chi(T_V(\langle D \rangle)\leq 0$. Then $d=3k-5$ where $k\geq 5$, $C$ has exactly three cusps with multiplicity sequences  $[(3(k-3),3_{k-3}),(4_{k-3},2_2),(2)]$ and $\chi(T_V(\langle D \rangle)= 0$. 
Moreover,
up to projective equivalence, this curve $C=C_k$  is given by a parametrization
$$(s:t) \mapsto (t^{3k-6}(t-s):s^8\tilde q_k(s,t):s^4t^{3k-9}),$$
where $\tilde q_k \in \C[s,t]$ is the homogeneous polynomial of degree $3k-13$
defined below.
\end{thm}
For the description of the polynomial $\tilde q_k$ we use the fomulas given in \cite[Main Theorem (h)]{BZ07}, following the suggestion kindly offered to us by Karol Palka. Define first a sequence of polynomials $X_k \in \C[u]$ by the formulas $X_4(u)=2u^2-u^3$ and
$$X_{k+1}(u)=\frac{u^4(X_k(u)-X_k(1))}{u-1}.$$
Note that $\deg X_k=3(k-3)$ and $X_k$ is divisible by $u^4$. Then set
$$Y_k(t)=t^{3(k-3)}X_k(\frac{1}{t}),$$
denote by $\tilde Y_k(s,t)$ the polynomial obtained by homogenization of 
the polynomial $Y_k$, and note that $\tilde Y_k(s,t)$ is divisible by $s^4$.
Set $\tilde q_k=\tilde Y_k(s,t)/s^4$. For instance, we get in this way the following parametrizations.

\begin{equation}\label{eqK5}
\tilde q_5(s,t)=-s^2+st+t^2.
\end{equation}

\begin{equation}\label{eqK6}
\tilde q_6(s,t)=-s^5+s^3t^2+s^2t^3+st^4+t^5.
\end{equation}

\begin{equation}\label{eqK7}
\tilde q_7(s,t)=-s^8-s^7t+s^5t^3+2s^4t^4+3s^3t^5+3s^2t^6+3st^7+3t^8.
\end{equation}

\begin{conj}\label{conj3C}
The curves $C_{k}\in \F\E(k)$  are free divisors and $mdr(f)=k-1$ for all  $k \geq 5$.
In particular, $\tau(C_k)=7(k-2)^2-k+3$  for any $k\geq 5$.
\end{conj}
This conjecture was verified by direct computation for $5 \leq k \leq 7$ using the above parametrizations.

The above conjectures can be summarized in the following main conjecture.

\begin{conj}\label{conjMC}
Any rational cuspidal plane curve $C$ of degree $d \geq 6$ and having at least three cusps is free.
Moreover, its exponents $(d_1,d_2)$ determine the infinite series to which it belongs, namely
\begin{enumerate}

\item $C\in  \F\ZZ_1(d,a)$ if and only if $d_1=2$, and in this case $d_2=d-3$.

\item $C\in  \F\ZZ_2(k)$ if and only if $d_1=(d-1)/2 \geq 3$, and in this case $d_2=d_1$. 

\item $C\in  \F\E(k)$ if and only if $d_1=(d+2)/3 \geq 4$, and in this case $d_2=(2d-5)/3$.

 \end{enumerate}

\end{conj}

Note that $d_1=mdr(f)$, so the above conjecture can be restated using $mdr(f)$. In view of formula \eqref{eqTAU1}, the above conjecture can be restated using the total Tjurina number $\tau(C)$ as well.

\subsection{On the rigidity conjecture}

We recall here the rigidity conjecture, following \cite{FZ96}. Let $C \subset \PP^2$ be a reduced, irreducible curve and let $V \to \PP^2$ be the minimal embedded resolution of the singularities of $C$, such that the total transform $D$ of $C$ is a SNC-divisor in $V$. Let $T_V\langle D \rangle$ be the logarithmic tangent bundle of the pair $(V,D)$. Then 
$H^0(T_V\langle D \rangle)$ is the space of infinitesimal automorphisms of the pair $(V,D)$, $H^1(T_V\langle D \rangle)$ is the space of the infinitesimal deformations  of the pair $(V,D)$ and $H^2(T_V\langle D \rangle)$ is the space of the obstructions for extending such infinitesimal deformations. Moreover the deformations of the pair $(V,D)$ corresponds to equisingular deformations of $C$ inside $\PP^2$. We say that $C$ is FZ-projectively rigid if $h^1(T_V\langle D \rangle)=0$ and that is
FZ-unobstructed if $h^2(T_V\langle D \rangle)=0$. Flenner and Zaidenberg show in \cite{FZ96}, see bottom of page 444, that in fact $C$ is FZ-projectively rigid if and only if $C$ is $es$-rigid as in our Definition \ref{def1}. They continue by stating the following conjecture.

\begin{conj}[Rigidity Conjecture]
\label{c2}
A rational cuspidal curve $C \subset \PP^2$ which has at least 3 singularities is FZ-projectively rigid and
FZ-unobstructed. In particular, such a curve is $es$-rigid.
\end{conj}
Flenner and Zaidenberg have checked in \cite{FZ96,FZ00} that this conjecture holds for the series of cuspidal rational plane curves $\F\ZZ_1(d,a)$ and $\F\ZZ_2(k)$. One has also the following conjecture, also due to Flenner and Zaidenberg.
\begin{conj}[Weak Rigidity Conjecture]
\label{c2w}
A rational cuspidal curve $C \subset \PP^2$ which has at least 3 singularities satisfies $\chi(T_V\langle D \rangle)=0$.
\end{conj}
Fenske has shown in \cite{Fe99} that the curves in the series $\F\E(k)$ satisfy the weak rigidity conjecture. A stronger conjecture, called the Negativity Conjecture, was proposed by Palka \cite{Pa1,Pa2}, and has extremely interesting consequences for the study of rational cuspidal plane curves.

On the other hand, our Conjecture \ref{conjMC} implies that 
any rational cuspidal plane curve $C$ of degree $d \geq 6$ and having at least three cusps is $ea$-rigid, by Corollary \ref{cor1} (1), but most of them are not $ea$-unobstructed by Corollary \ref{cor3} (1).

The results by Flenner, Zaidenberg and Fenske mentioned above imply that for a curve $C$ in any of the families $ \F\ZZ_1(d,a)$,  $  \F\ZZ_2(k)$ and   $\F\E(k)$, one has the equality
$$G \cdot C=\V_d^{ea}(S_1,S_2,S_3)=\V_d^{es}(S_1,S_2,S_3),$$ because any equisingular deformation of $C$ is in fact an equianalytic deformation. This equality
holds also for all the curves in Proposition \ref{prop0C}.
Note that the equality $\V_d^{ea}(S_1,...,S_r)=\V_d^{es}(S_1,...,S_r)$ fails for rational cuspidal curves with one cusps, as Examples \ref{ex35} and \ref{ex36} below shows. In fact, the rational unicuspidal curves constructed in Theorem 1.1 and Theorem 1.2 in \cite{DStExpo} belong, for a fixed degree $d$, to a unique analytic subspace $\V_d^{es}(S_1)$, but to distinct analytic subspaces $\V_d^{ea}(S_1)$, as their Tjurina numbers vary.
 Similar examples for curves with two cusps are also easy to produce,
for instance starting with the curve $C: f=y^7+x^3z^4=0$.
On the other hand, the equality $G \cdot C=\V_d^{ea}(S_1,...,S_r)$
holds for rational cuspidal curves of degree $d \geq 6$ which satisfy in addition the condition $\mu(C)=\tau(C)$, see \cite{Drcc}.

\section{Jacobian syzygies, Jacobian module and equisingular deformations} 

Recall the ideal sheaves $\I(Z^*(C))$ introduced before Lemma \ref{lem1}. They enter into the following exact sequence
\begin{equation}\label{ES10}
0 \to \I(Z^{ea}(C)|C)\to \I(Z^{es}(C)|C) \to \F_C \to 0,
\end{equation}
where the sheaf $\F$ is just the quotient o sheaf $\I(Z^{es}(C)|C)/\I(Z^{ea}(C)|C)$. It follows that the support of the sheaf $\F_C$ is contained in the set of non-simple singularities of $C$, and at each such singularity $p$ one has
\begin{equation}\label{ES11}
\dim  \F_{C,p} = \dim I^{es}(C,p)/I^{ea}(C,p)=m^{es}(C,p),
\end{equation} 
as in \eqref{eq01}. To get the value of this dimension in concrete cases, one may use the algorithm described in \cite{CGL} and implemented in SINGULAR \cite{Sing}, or Example \ref{ex1} in very simple cases. As in Lemma \ref{lem1}, we have
$\I(Z^{ea}(C)|C)=\J_f|C$, and hence we get the following long exact sequence of cohomology groups
\begin{equation}\label{ES12}
0 \to H^0(\J_f|C)\to H^0(\I(Z^{es}(C)|C)) \to H^0(\F_C) \to H^1(\J_f|C)\to H^1(\I(Z^{es}(C)|C)) \to 0.
\end{equation}
This yields the following result.
\begin{thm}
\label{thm10}
Let $C:f=0$ be a reduced plane curves such that $C \in V_d^{ea}(S_1,...,S_r)$. Then  the Zariski tangent space
 $T_C \V_d^{es}(S_1,...,S_r)$ of the analytic subspace
$ \V_d^{es}(S_1,...,S_r)$ at the point $C$ sits in the following exact sequence
$$0  \to T_C \V_d^{ea}(S_1,...,S_r) \to  T_C \V_d^{es}(S_1,...,S_r) \to H^0(\F_C).$$
Moreover, the morphism $T_C \V_d^{es}(S_1,...,S_r) \to H^0(\F_C)$ is surjective if  $mdr(f)\geq d-4$.
In particular, one has
$$\dim T_C \V_d^{ea}(S_1,...,S_r) \leq \dim T_C \V_d^{es}(S_1,...,S_r)  \leq \dim T_C \V_d^{ea}(S_1,...,S_r)+ \sum_p m^{es}(C,p),$$
and equality holds in the second inequality if  $mdr(f)\geq d-4$. Moreover, the analytic subspace
$ \V_d^{es}(S_1,...,S_r)$ is smooth of the expected dimension at the point $C$ if  $mdr(f)\geq d-4$.
\end{thm}
Here the sum $\sum_p m^{es}(C,p)$ is over all the non-simple singularities of $C$. The last claim in this theorem follows from the formula \eqref{eq0}, Theorem \ref{thm1} (2) and the identification $H^1(\J_f|C)=AR(f)_{d-5}$ in the exact sequence \eqref{ES12}, recall the proof of Theorem \ref{thm4}.

\begin{ex}
\label{ex35}
Consider the rational cuspidal curve
$$C: f=y^5+x^4z=0.$$
This curve has a unique singularity $S_1$ at $p_1=(0:0:1)$ with $\mu(C,p_1)=12$.
The Jacobian ideal $J_f$ is generated by $x^4,y^4,x^3z$ and the
ideal $I_f$ is generated by $x^3,y^4$. It follows that $h=x^3y^2 \in I_f$, but $h \notin J_f$, hence $h$ yields a basis for $N(f)_5$. For $t \ne 0$, the curve 
$$C_t:f_t=f+th=y^5+x^4z+tx^3y^2=0$$
has again a singularity at $p_1=(0:0:1)$, which is analytically isomorphic to the singularity $S_1$. It follows that the deformation $(C_t,p_1)$ is equianalytic. If we apply Theorem \ref{thm4}, we see that the analytic subspace
$\V_5^{ea}(S_1)$ is smooth of the expected dimension $ed$ at $C$, where
$ed=20-12=8$. It is easy to see that all the curves $C_t$ for $t \ne 0$ are in the same $G$-orbit. Moreover $\dim G \cdot C=7$, $\dim G \cdot C_1=8$ and $G \cdot C$ is obviously contained in the closure of the orbit $G \cdot C_1$. It follows that one has the equality
$$\V_5^{ea}(S_1)= G \cdot C \cup G \cdot C_1$$
as germs at $C$. Note that $C$ is a nearly free curve with $mdr(f)=1$, while $C_1: f_1=0$ is a free curve with $mdr(f_1)=2$.

\medskip

Now we turn our attention to the analytic subspace $\V_5^{es}(S_1)$. Using Example \ref{ex1}, we see that the equisingular ideal $I^{es}$ of the singularity
$T_{4,5}:u^4+v^5=0$ is spanned by $u^3,v^4$ and $u^2v^3$, in particular
$m^{es}(C,p)=1$. On the other hand, it is known that the tangent space $T_C(G \cdot C)$ is given by the projectivization of $J_{f,5}$, and hence a transversal $\T$ to this tangent space is spanned by the following 
$${7 \choose 2}-\dim J_{f,5}=21-8=13$$
monomials: $x^3y^2$, $x^2y^az^b$ for $a,b \geq 0$, $a+b=3$,
$xy^iz^j$ for $i \geq 0$, $j \geq 1$, $i+j=4$, $y^pz^q$, $p\geq 0$, $q \geq 2$, $p+q=5$. Any deformation of $C:f=0$ is given by a deformation in this transversal $\T$, up-to the action of the group $G$.
If we transform the above monomials in $x,y,z$ into monomials in $u,v$ by setting $x=u$, $y=v$ and $z=1$, and we look for those monomials in $u,v$ giving $\mu$-constant deformations of the singularity $T_{4,5}$, we get only the following two monomials $x^3y^2$ and $x^2y^3$.
The deformation in the direction of the monomial $x^3y^2$ was already discussed above. For $t \ne 0$, the curve 
$$C'_t:f'_t=y^5+x^4z+tx^2y^3=0$$
has again a singularity at $p_1=(0:0:1)$, which is equisingular to the singularity $S_1$, since the deformation is $\mu$-constant. On the other hand, this deformation is no longer analytically trivial, since 
$$\tau(C'_1,p_1)=11<12=\tau(C,p_1)=\mu(C,p_1).$$
Indeed, all the curves $C'_t$ for $t \ne 0$ are in the same $G$-orbit, and one can check that $C_1'$ is nearly free with $mdr(f)=2$. 

If we apply Theorem \ref{thm10}, we see that the analytic subspace
$\V_5^{es}(S_1)$ is smooth of the expected dimension $ed$ at $C$, where
$ed=20-12+1=9$.
It follows that we have the following equality of smooth germs at $C$
$$\V_5^{es}(S_1)=G \cdot (f+\T_0),$$
where $\T_0$ is the 2-dimensional vector subspace in the transversal $\T$ spanned by the two monomials $x^3y^2$ and $x^2y^3$. In fact, note that the curve
$$C_{\al, \be}:f+\al x^3y^2+ \be x^2y^3=0$$
is projectively equivalent to the curve $C_{\al t^3, \be t^2}$, for any $t \ne 0$. It follows that
$$\V_5^{es}(S_1)=G \cdot C \cup_{[\al:\be] \in \PP(3,2)}G \cdot C_{\al, \be},$$
i.e. $\V_5^{es}(S_1)$ is the union of the 7-dimensional orbit $G \cdot C$ and of a 1-parameter family of 8-dimensional $G$-orbits, all of them containing $G \cdot C$ in their closure, parametrized by the weighted projective line $\PP(3,2)=\PP^1$.

\end{ex}

\begin{rk}
\label{rk35}
Note that the deformation $C'_t$ above has essentially the same properties as the deformation 
$$C''_t:f''_t=y^5+x^4z+txy^4=0$$
in the direction of $xy^4 $, a vector in the tangent space $T_C(G \cdot C)$.

\end{rk}

\begin{ex}
\label{ex36}
Consider the rational cuspidal curve
$$C: f=y^6+x^5z=0.$$
This curve has a unique singularity $S_1$ at $p_1=(0:0:1)$ with $\mu(C,p_1)=20$.
The Jacobian ideal $J_f$ is generated by $x^5,y^5,x^4z$ and the
ideal $I_f$ is generated by $x^4,y^5$. It follows that $h=x^4y^2 \in I_f$, but $h \notin J_f$, hence $h$ yields a basis for $N(f)_6$. For $t \ne 0$, the curve 
$$C_t:f_t=f+th=y^5+x^4z+tx^4y^2=0$$
has again a singularity at $p_1=(0:0:1)$, which is equisingular to the singularity $S_1$, but has 
$$\tau(C_t,p_1)=19<20=\tau(C,p_1)=\mu(C,p_1).$$ 
Hence a deformation given by a monomial in $I_f$ is not necessarily equianalytic as in Example \ref{ex35} above.  Note that all the curves $C_t$ for $t \ne 0$ are in the same $G$-orbit,
$C$ is a nearly free curve with $mdr(f)=1$, while $C_1: f_1=0$ is a free curve with $mdr(f_1)=2$.

 If we apply Theorem \ref{thm3}, we see that 
$\dim T_C\V_6^{ea}(S_1)=8$. Moreover, we know by \cite[Theorem 1.1]{Drcc} that one has $V_6^{ea}(S_1)=G \cdot C$, and hence 
$$\dim V_6^{ea}(S_1)=\dim G \cdot C=7.$$
Hence the analytic subspace $\V_6^{ea}(S_1)$ is non-reduced at any point, even though its support $V_6^{ea}(S_1)$ is smooth.
Note that in this example, unlike the classical example of B. Segre
of curves with many simple cusps, discussed in \cite{Ta} or in \cite[Example 4.7.10]{SeBook}, the dimension of $V_6^{ea}(S_1)$ is the expected dimension
$$d(d+3)/2-\tau(C)=27-20=7.$$

\medskip

Now we turn our attention to the analytic subspace $\V_6^{es}(S_1)$. Using Example \ref{ex1}, we see that the equisingular ideal $I^{es}$ of the singularity
$T_{5,6}:u^5+v^6=0$ is spanned by $u^4,v^5, u^2v^4$ and $u^3v^3$, in particular
$m^{es}(C,p)=3$. On the other hand, it is known that the tangent space $T_C(G \cdot C)$ is given by the projectivization of $J_{f,6}$, and hence a transversal $\T$ to this tangent space is spanned by the following 
$${8 \choose 2}-\dim J_{f,5}=28-8=20$$
monomials: $x^4y^2$, $x^3y^az^b$ for $a,b \geq 0$, $a+b=3$,
$x^2y^iz^j$ for $i,j \geq 0$,  $i+j=4$, $xy^pz^q$, $p\geq 0$, $q \geq 1$, $p+q=5$, and $y^mz^n$ for $m\geq 0$, $n \geq 2$, $m+n=6$. As above, any deformation of $C:f=0$ is given by a deformation in this transversal $\T$, up-to the action of the group $G$.
If we transform the above monomials in $x,y,z$ into monomials in $u,v$ by setting $x=u$, $y=v$ and $z=1$, and we look for those monomials in $u,v$ giving $\mu$-constant deformations of the singularity $T_{5,6}$, we get only the following three monomials $x^4y^2$, $x^3y^3$ and $x^2y^4$.

If we apply Theorem \ref{thm10}, we see that the analytic subspace
$\V_6^{es}(S_1)$ satisfies
$$8 \leq \dim T_C\V_6^{es}(S_1)<8+3=11.$$
Note that the curve
$$C_{\al, \be, \gamma}:f+\al x^4y^2+ \be x^3y^3+\gamma x^2y^4=0$$
is projectively equivalent to the curve $C_{\al t^4, \be t^3, \gamma t^2}$, for any $t \ne 0$. It follows that
$$\V_6^{es}(S_1)=G \cdot C \cup_{[\al:\be:\gamma] \in \PP(4,3,2)}G \cdot C_{\al, \be, \gamma},$$
i.e. $V_6^{es}(S_1)$ is the union of the 7-dimensional orbit $G \cdot C$ and of a 2-parameter family of 8-dimensional $G$-orbits, all of them containing $G \cdot C$ in their closure, parametrized by the weighted projective plane $\PP(4,3,2)$. It follows that 
$\dim V_6^{es}(S_1)=8+2=10$ and hence we get
 $$\dim T_C\V_6^{es}(S_1)\leq 10=\dim \V_6^{es}(S_1).$$
 It follows that the analytic subspace $\V_6^{es}(S_1)$ is smooth at $C$, and of the expected dimension, which is $ed=27-20+3=10$, although $1=mdr(f) <d-4=2$.
 
 Note that one has an isomorphism of weighted projective spaces $\PP(4,3,2)=\PP(2,3,1)$, see \cite{Del, Dol}, and this implies that the parameter space $X=\PP(2,3,1)$ has two singular points, namely
 $x_1=(1:0:0)$ and $x_2=(0:1:0)$, see \cite{DiDi}.
The curve corresponding to $x_1$ is the curve $C_1$ introduced above. Using Theorem \ref{thm10}, it follows that  the analytic subspace
$\V_6^{es}(S_1)$ is smooth of the expected dimension $ed_1$ at $C_1$, where
$ed_1=27-19+2=10$, since $m^{es}(C_1,p)=2$. The curve corresponding to $x_2$ is the curve $$C'_1:f_1'=f+x^3y^3=0,$$
which is nearly free with $mdr(f_1')=3$ and $\tau(C'_1,p_1)=18$.
Using Theorem \ref{thm10}, it follows that  the analytic subspace
$\V_6^{es}(S_1)$ is smooth of the expected dimension $ed'$ at $C'_1$, where
$ed'=27-18+1=10$, since $m^{es}(C'_1,p)=1$.
In this example the Tjurina numbers $\tau(C,p)$ and $\tau^{es}(C,p)$ and the Tjurina modalities $m^{es}(C,p)$ have been computed using the SINGULAR software, \cite{Sing}.
\end{ex}

\begin{rk}
\label{rk36}
Note that the deformation of the sextic $C$ above given by
$$C''_t:f''_t=y^6+x^5z+txy^5=0$$
in the direction of $xy^5 $, a vector in the tangent space $T_C(G \cdot C)$ produces nearly free curves with $mdr(f''_t)=2$ for $t \ne 0$.
One can show by simple linear transformations that this deformation is $G$-equivalent to a equisingular deformation of $C$ inside one of the orbits $G \cdot C_{\al, \be, \gamma}$.

\end{rk}

\end{document}